\documentclass{amsart}
\usepackage[english]{babel}
\usepackage{epsfig}
\usepackage{amssymb}
\usepackage{amsmath}
\usepackage{amsthm,youngtab}
\usepackage{enumerate}
\usepackage{cite}
\usepackage{xcolor}
\usepackage{url}

\title[Polynomials and the exponent of matrix multiplication]
{Polynomials and the exponent of \\
matrix  multiplication}
\date{}
\newcommand{\C}{\mathbb{C}}
\newcommand{\Z}{\mathbb{Z}}
\newcommand{\PP}{\mathbb{P}}

\newcommand{\Qc}{\mathcal{Q}}

\newcommand{\sL}{\mathcal{L}}

\newcommand{\rank}{\operatorname{rank}}

\newtheorem{defn0}{Definition}[section]
\newtheorem{conj}[defn0]{Conjecture}
\newtheorem{prop0}[defn0]{Proposition}
\newtheorem{prop0*}[defn0]{Proposition*}
\newtheorem{thm0}[defn0]{Theorem}
\newtheorem{thm0*}[defn0]{Theorem*}
\newtheorem{lemma0}[defn0]{Lemma}

\newtheorem{rem0}[defn0]{Remark}

\newcommand{\tdim}{\operatorname{dim}}
\newcommand{\dieci}{{10}}

\def\Mn{M_{\langle n \rangle}}\def\Mtwo{M_{\langle 2\rangle}}
\def\SMn{sM_{\langle n \rangle}}\def\LMn{  \Lambda M_{\langle n \rangle}}\def\SMthree{sM_{\langle 3\rangle}}
\def\SMsn{sM_{\langle n \rangle}^S}\def\SMson{sM_{\langle n \rangle}^{S,0}}
\def\ZMn{sM_{\langle n \rangle}^{Z}}
\def\SMtwo{sM_{\langle 2 \rangle}}
\def\sSMtwo{sM_{\langle 2 \rangle}^S}

\def\Matn{{\rm{Mat}}_n} 
\def\ttrace{\operatorname{trace}}
\def\ot{\otimes}
\def\BZ{\Bbb Z}\def\BC{\mathbb C}
\def\La#1{\Lambda^{#1}}
\def\ur{\underline {\bold R}}
\def\tend{\operatorname{End}}
\def\op{{\mathord{\,\oplus }\,}}\def\o{\omega}
\def\tlog{\operatorname{log}}
\def\fgl{\mathfrak g\mathfrak l}
\def\trank{\operatorname{rank}}\def\ttrace{\operatorname{trace}}
\def\s{\sigma}\def\pp#1{\mathbb P^{#1}}\def\FS{\mathfrak  S}
\newcommand{\BP}{\mathbb{P}}

\author[L.~Chiantini]{Luca Chiantini}
\address{Dipartimento di Ingegneria dell'Informazione e Scienze Matematiche, 
Universit\`a di Siena, Italy}
\email{luca.chiantini@unisi.it}

\author[J.~Hauenstein]{Jonathan D. Hauenstein}
\address{Department of Applied and Computational Mathematics and Statistics, University of Notre Dame, Notre Dame, IN, USA}
\email{hauenstein@nd.edu}

\author[C.~Ikenmeyer]{Christian Ikenmeyer}
\address{Max Planck Institute for Informatics, Saarland Informatics Campus, Germany}
 \email{cikenmey@mpi-inf.mpg.de}

\author[J.M.~Landsberg]{J.M. Landsberg}
\address{Department of Mathematics, Texas A\&M University, College Station, TX, USA}
 \email{jml@math.tamu.edu} 

\author[G.~Ottaviani]{Giorgio Ottaviani}
\address{Dipartimento di Matematica e Informatica ``Ulisse Dini'', Universit\`a di 
Firenze, Italy}
\email{ottavian@math.unifi.it}

\numberwithin{equation}{section}    
\begin{document}

\begin{abstract}
We define tensors, corresponding to cubic polynomials, which have the same 
exponent $\omega$ as the matrix multiplication tensor. 
In particular, we study
the symmetrized  matrix multiplication tensor $\SMn$
defined on an $n\times n$ matrix $A$ by $\SMn(A)=\ttrace(A^3)$.
The use of polynomials enables the introduction of 
additional techniques from algebraic geometry in the study 
of the matrix multiplication exponent $\omega$. 
\end{abstract}

\maketitle

\section{Introduction}

The {\it exponent of matrix multiplication} is the smallest constant $\o$
such that two
  $n\times n$ matrices   may be multiplied 
  by performing   
$O(n^{\omega+\epsilon})$
  arithmetic operations for every $\epsilon>0$. 
  It is a central
   open problem to estimate $\omega$ since it
governs the complexity of many basic algorithms in linear algebra.
The current state of the art \cite{copwin135,williams,legall,stothers}~is 
$$2\le\omega< {2.374}.$$

A  {\it tensor} $t\in \BC^N\otimes \BC^N\otimes \BC^N$ has {\it (tensor) rank} $r$ if
$r$ is the minimum such that there
exists $u_i,v_i,w_i\in\BC^N$ with
$t=\sum_{i=1}^ru_i\ot v_i\ot w_i$. 
In this case, we write ${\bold R}(t)=r$.
Let $V=\BC^n$ and $\tend (V)=\Matn$ be the vector space of $n\times n$ matrices over $\C$.
The matrix multiplication tensor $\Mn\in \Matn^\vee\otimes \Matn^\vee\otimes \Matn^\vee$
is
\begin{equation}\label{eq:deftn}\Mn(A,B,C)=\ttrace(ABC),\end{equation}
{where $\Matn^\vee$ is the vector space dual to $\Matn$.}

Strassen \cite{Strassen505} showed  that $\omega=\liminf[\tlog_n(\bold R(\Mn))]$.
If the tensor $t$ can be expressed as a limit of tensors of 
rank $s$ (but not a limit of tensors of rank at 
most~\mbox{$s-1$}), then~$t$ has {\it border rank} $s$, denoted $\ur(t)=s$.
This is equivalent to $t$ being in the Zariski closure of the 
set of tensors of rank $s$ but not in the Zariski closure of the set 
of tensors of rank at most $s-1$, see, e.g., \cite[Thm. 2.33]{MR1344216}.
This was rediscovered in complexity theory  in \cite{ald:84}. Bini \cite{MR605920} showed that
$\omega=\liminf[\tlog_n(\ur(\Mn))]$.

The determination of the fundamental constant $\o$ is a central question in algebraic complexity theory.
In 1981, Sch\"onhage~\cite{MR623057} 
showed the exponent $\omega$ could
be bounded using disjoint sums of matrix multiplication
tensors.  Then, in 1987, Strassen \cite{MR882307} proposed using tensors other than $\Mn$ which are 
easier to analyze
due to their combinatorial properties to prove upper bounds on $\o$.
These other tensors are then degenerated
to disjoint matrix multiplication tensors. 
The main goal of this paper is to open a different path to bounding $\omega$
by introducing polynomials that are closely related to 
matrix multiplication.

We expect these polynomials are easier to work with
in two ways.  First, we want to take advantage of the vast literature in algebraic geometry regarding the
geometry of cubic hypersurfaces.  
Second, we want to exploit recent numerical computational techniques. 
The difficulty of the usual matrix multiplication
tensor is the sheer size of the problem,
even for relatively small $n$.
Despite considerable effort, no $4\times 4$ decompositions, other than the standard
rank 64 decomposition
and the rank 49 decomposition obtained by squaring Strassen's $2\times 2$ decomposition, have
appeared in the literature.
With our approach, the polynomials are defined
on much smaller spaces thereby allowing
one to perform more computational experiments 
and produce additional data for forming conjectures.

Let $\mathrm{Sym}^3\BC^N\subset (\BC^N)^{\ot 3}$ 
and $\La 3\BC^N\subset (\BC^N)^{\ot 3}$ 
respectively denote the space of symmetric 
and skew-symmetric tensors.
{Tensors in $\mathrm{Sym}^3\BC^N$ 
may be viewed as homogeneous  cubic polynomials 
 in $N$ variables.}
While the matrix multiplication tensor $\Mn$ 
is neither symmetric nor skew-symmetric,
it is $\BZ_3$-invariant 
where $\BZ_3$ denotes the cyclic group on three elements permuting the factors
since $\ttrace(ABC)=\ttrace(BCA)$.
The space of $\BZ_3$-invariant tensors in~$(\BC^N)^{\ot 3}$ is
$$
[(\BC^N)^{\ot 3}]^{\BZ_3}= \mathrm{Sym}^3\BC^N  \op \La 3\BC^N.
$$
Thus, respectively define the 
symmetrized and skew-symmetrized part
of the matrix multiplication tensor, namely
\begin{align}
\label{eq:smn} \SMn(A,B,C)&:=\frac 12[ \ttrace(ABC)+\ttrace(BAC)] \\
 \LMn(A,B,C)&:= \frac 12[\ttrace(ABC)-\ttrace(BAC) ] \end{align}
so that 
\begin{equation}
\label{msml}
\Mn=\SMn+\LMn.
\end{equation}
The $\BZ_3$-invariance implies
$\SMn\in \mathrm{Sym}^3\BC^N$ and
$\LMn\in \La 3\BC^N$.

The tensor $\Mn$ is the structure tensor for the algebra $\Matn$. {Similarly}, the skew-symmetrized matrix multiplication tensor $\LMn$ is (if one ignores
the $\frac 12$) the structure
tensor for the Lie algebra $\fgl(V)$.
The symmetrized matrix multiplication
tensor $\SMn$ 
is the structure tensor for $\Matn$ considered as a Jordan algebra, i.e., with the multiplication $A\circ B=\frac 12(AB+BA)$.  In particular, considered as a cubic polynomial on $\Matn$, $$\SMn(A)=\ttrace(A^3).$$

We further define the following cubic polynomials (symmetric tensors):
\begin{itemize}
\item $\SMsn$: restriction of $\SMn$ to symmetric matrices $\mathrm{Sym}^2V$, 
\item $\SMson$: restriction of $\SMsn$ to traceless symmetric matrices, and 
\item $\ZMn$: restriction of $\SMsn$ to symmetric matrices with
zeros on diagonal.
\end{itemize}

In order to have an invariant definition of $\SMsn$ and $\SMson$, one needs an identification of $V$ with $V^*$.
{ Two natural ways of obtaining this identification are via a
nondegenerate symmetric quadratic form or, when $\tdim V$ is even, a
skew-symmetric form. We will often use the former, which reduces the symmetry group from the general
linear group to the orthogonal group.}
We do not know of a nice invariant definition for the polynomial $\ZMn$.

For a homogeneous degree $d$ polynomial
  $P$, the {\it symmetric or Waring rank} $\bold R_s(P)$  is  the smallest
$r$ such that $P=\sum_{j=1}^r \ell_j^d$, where $\ell_j$ are linear forms.  The {\it symmetric border rank}  $\ur_s(P)$ is  the smallest $r$ such that $P$ is
a limit of polynomials of {symmetric} rank at most $r$.
Note that
\begin{equation}\label{obviousinq} \bold R(P)\leq \bold R_s(P) \hbox{~~~and~~~}\ur(P)\leq \ur_s(P).
\end{equation} 
We notice that there are several general cases where equality holds
in both of these relations. We refer to \cite{Comon2008,MR3092255} for a discussion.

Our main result is that one can compute the exponent $\o$
of matrix multiplication using these polynomials 
even when considering symmetric rank and border rank.

\begin{thm0}\label{thm:main}Let $\o$ denote the exponent of matrix multiplication. Then
\begin{align}
\omega & 
  \label{SMnom}=\liminf_n \left[\log_n {\bold R}(\SMn)\right]=\liminf_n \left[\log_n \ur(\SMn)\right]
  \\
  \nonumber &= 
   \liminf_n \left[\log_n {\bold R_s}(\SMn)\right]  
 = \liminf_n \left[\log_n \ur_s(\SMn)\right] \\
   &
   \label{SMsnom} =\liminf_n \left[\log_n {\bold R}(\SMsn)\right]=\liminf_n \left[\log_n \ur(\SMsn)\right]
    \\
  \nonumber &=  
   \liminf_n \left[\log_n {\bold R_s}(\SMsn)\right]  
 = \liminf_n \left[\log_n \ur_s(\SMsn)\right] \\  
   &\label{SMsoom}=\liminf_n \left[\log_n {\bold R}(\SMson)\right]=\liminf_n \left[\log_n \ur(\SMson)\right]
    \\
  \nonumber &=  
   \liminf_n \left[\log_n {\bold R_s}(\SMson)\right]  
 = \liminf_n \left[\log_n \ur_s(\SMson)\right] \\  
    &\label{ZMnom} =\liminf_n \left[\log_n {\bold R}(\ZMn)\right]=\liminf_n \left[\log_n \ur(\ZMn)\right]
     \\
  \nonumber &=  
   \liminf_n \left[\log_n {\bold R_s}(\ZMn)\right]  
 = \liminf_n \left[\log_n \ur_s(\ZMn)\right].
\end{align}
\end{thm0}

Proofs are   given in \S\ref{sec:symrank} for \eqref{SMnom}, 
\S\ref{sec:smsn} for (\ref{SMsnom}) and (\ref{SMsoom}), and \S\ref{sec:szmn}
for \eqref{ZMnom}.

\subsection{Explicit ranks and border ranks}

For any $t\in \BC^N\ot \BC^N\ot \BC^N$, 
  the symmetrization of $t$ 
is $\mathcal{S}(t):=\frac 1 6 \sum_{\pi\in\mathfrak{S}_3}\pi(t)\in \mathrm{Sym}^3\BC^N$.
In particular, 
$\mathcal{S}(t)=t$ 
if and only if $t\in \mathrm{Sym}^3\BC^N$.
The following provides bounds relating $t$ and $\mathcal{S}(t)$.

\begin{lemma0}\label{lem:symrk} For $t\in \BC^N\ot \BC^N\ot \BC^N$,  
$ {\bold R}_s(\mathcal{S}(t))\le 4{{\bold R}}(t)$ and 
$ \ur_s(\mathcal{S}(t))\le 4{\ur}(t)$.
\end{lemma0}
\begin{proof}  If 
$t=\sum_{i=1}^ru_i\ot v_i\ot w_i$ with $u_i, v_i, w_i, \in V$,
then $\mathcal{S}(t)=\sum_{i=1}^r (u_i  v_i  w_i)$. Since  ${\bold R}_s(xyz)=4$
(see, e.g., \cite[\S 10.4]{Land_book}),
this immediately yields that ${\bold R}_s(\mathcal{S}(t))\le 4{{\bold R}}(t)$. In the same way, if $t$ is a limit of tensors of the form $\sum_{i=1}^ru_i\ot v_i\ot w_i$, 
this yields $ \ur_s(\mathcal{S}(t))\le 4{\ur}(t)$.
 \end{proof}

In particular, $\bold R(\SMn)\leq 2\bold R(\Mn)< 2n^3$ (as $\SMn$ is the sum of two matrix multiplications, by (\ref{eq:smn})) so that 
$\bold R_s(\SMn)\leq 8\bold R(\Mn)< 8n^3$ and similarly for all its degenerations. 

The following summarizes some results about small cases.

\begin{thm0}\label{smallcasesrk}
\mbox{ }
\begin{enumerate}
\item\label{item:1}
 $ {\bold R}_s(\SMtwo)=6$ and $\ur_s(\SMtwo)=5$
  (\hspace{1sp}\cite[IV, \S 97]{Seg}, \cite[Prop. 7.2]{LT}).
\item\label{item:3}  
$\ur_s(sM_{\langle 3\rangle})\ge 14$.
\item\label{item:2}
${\bold R}_s(\sSMtwo)=\ur_s(\sSMtwo) = 4$
(\hspace{1sp}\cite[IV \S 96]{Seg} or \cite[\S 8]{LT}).
\item\label{item:4}  
$\ur_s(sM^{S}_{\langle 3\rangle}) = 10$.
\item\label{item:5} $sM^{S,0}_{\langle 2\rangle} = 0$ while $\ur_s(sM^{S,0}_{\langle 3\rangle})= 
\bold R_s(sM^{S,0}_{\langle 3\rangle}) = 8$.
\item\label{item:6} $\ur_s(sM^{S,0}_{\langle 4\rangle})\ge 14$.
\item\label{item:7} $sM^{Z}_{\langle 2\rangle} = 0$
while 
$\ur_s(\ZMn) =  \bold R_s(\ZMn)=2^{n-1}$ for $n=3,4,5$, 
with $\bold R_s(sM^Z_{\langle 6\rangle})\le 30$, 
$\bold R_s(sM^Z_{\langle 7\rangle})\le 48$, and $\bold R_s(sM^{Z}_{\langle 8\rangle})\le 64$.
\end{enumerate}
 \end{thm0} 
 
The cases (\ref{item:1}) and (\ref{item:2}) are discussed respectively in \S\ref{subsec:SMtwo} and \S\ref{subsec:smstwo}.
The case (\ref{item:4}) is proved in \S\ref{subsec:tableau} with a tableau evaluation. The cases (\ref{item:3}), (\ref{item:5}), (\ref{item:6})
 are proved with the technique of Young flattenings introduced in \cite{LOann}
which has already been 
used in the unsymmetric case in \cite{LO}. 
In particular, Proposition~\ref{prop:14} below 
considers (\ref{item:3}) with the other cases 
following analogously.
The case (\ref{item:7}) is proved by exhibiting explicit decompositions in Theorems \ref{thm:waringzmn1}, \ref{thm:waringzmn2}, and \ref{pro:waringzmn3}.

Since one of our goals is to simplify the problem
in order to further exploit numerical computations, we experiment with numerical tools and probabilistic
methods via {\tt Bertini} \cite{BertiniSoftware}.
We believe the computations could likely be 
converted to rigorous proofs, e.g., by showing
that an overdetermined system has a solution 
nearby the given numerical approximation \cite{Certification}.
We write \textbf{Theorem*} when we mean the result of a
numerical computation.

\begin{thm0*}\label{smallcasesrk*} 
$ {\bold R}_s(\SMthree)\leq 18$.
\end{thm0*} 
 
We show this in Theorem*~\ref{thm*:ms3_19}
with data regarding this and other 
computations available 
at \url{http://dx.doi.org/10.7274/R0VT1Q1J}.
 
\subsection*{Notation and conventions} The group of invertible linear maps $\BC^N\rightarrow \BC^N$ is
denoted $GL_N$ and the permutation group on $d$ elements by $\FS_d$. 
For $u,v,w\in\BC^N$, we
have $u\ot v\ot w\in (\BC^N)^{\ot 3}$ and
$uvw\in \mathrm{Sym}^3\BC^N$. 
The space $\Matn$ is canonically self-dual. Given a matrix $L$,
when we consider 
  $L\in \Matn^\vee$, 
we write $L^3\in \mathrm{Sym}^3(\Matn^\vee)$ 
for the cubic polynomial function which sends the matrix $A$~to~$\left[\mathrm{trace}(L^T A)\right]^3$,
where $L^T$ is the transpose of $L$. Note that $L^3$ is
a function and {\it not} the cube of the matrix~$L$.
In particular, $\SMn=\sum_{i=1}^k  L_i^3$ means that 
$$\mathrm{trace}(A^3)=\sum_{i=1}^k\left[\mathrm{trace}(L_i^T A)\right]^3.$$
For a partition $\pi$ of $d$, $S_{\pi}\BC^N$ denotes the corresponding $GL_N$-module and $[\pi]$ the corresponding $\FS_d$-module.
In particular $S_{(d)}\BC^N=\mathrm{Sym}^d\BC^N$ and $S_{(1^d)}\BC^N=\La d\BC^N$. 

\subsection*{Acknowledgement} This project began during the Fall 2014 program
{\it Algorithms and Complexity in Algebraic Geometry}. The authors thank
the Simons Institute for providing a wonderful research environment.

\section{The polynomial $\SMn$}
\label{sec:symrank}
 
We start with the first statement from Theorem~\ref{thm:main}.

\begin{proof}[Proof of \eqref{SMnom}]
  Lemma \ref{lem:symrk} and  \eqref{obviousinq}
imply 
$$ 4{\bold R}(\Mn)\ge{\bold R}_s(\SMn)\ge {\bold R}(\SMn)$$
so that
$$\omega\ge\liminf_n \left[\log_n {\bold R}_s(\SMn)\right] \ge \liminf_n \left[\log_n {\bold R}(\SMn)\right] .$$

For $n\times n$ matrices $A$, $B$, $C$ consider the $3n \times 3n$ matrix 
$X =\begin{pmatrix}0& 0& A\\
 C& 0& 0\\
 0& B& 0\end{pmatrix}$.
  Then, 
$X^3=\begin{pmatrix}ABC& 0& 0\\
 0& CAB& 0\\
 0& 0& BCA\end{pmatrix}$
and $\ttrace(X^3) = 3\ttrace(ABC)$.
This shows that
${\bold R}(\Mn)\le{\bold R}(sM_{\langle 3n\rangle})$
yielding the inequality 
$\omega\le \liminf_n \left[\log_n {\bold R}(\SMn)\right]$.
The border rank statement follows similarly by taking limits.
\end{proof}
 
As a $GL_N$-module via the Cauchy formula, $\mathrm{Sym}^3(\tend(V))=\mathrm{Sym}^3(V\ot V^*)$ decomposes as
\begin{align}
\mathrm{Sym}^3(\tend(V))&=\mathrm{Sym}^3V\ot \mathrm{Sym}^3V^*\op S_{21}V\ot S_{21}V^*\op \La 3 V\ot \La 3 V^*\\
&= \label{eq:sym3mn}\tend(\mathrm{Sym}^3V)\op \tend(S_{21}V)\op \tend(\La 3 V).
\end{align}
The tensor $\Mn\in \Matn^\vee\otimes \Matn^\vee\otimes \Matn^\vee=\mathrm{End}(V^{\otimes 3})$ corresponds to the identity endomorphism.
Since $V^{\otimes 3}=\mathrm{Sym}^3V\op  (S_{21}V)^{{ \op }2}\op \La 3 V$, it 
follows that $\mathrm{End}(V^{\otimes 3})$, as a $GL(V)$-module, contains the submodule
$$\tend(\mathrm{Sym}^3V)\op \left(\tend(S_{21}V)^2\right)\op \tend(\La 3 V).$$
The projection of $\SMn$ 
onto each of the three summands in (\ref{eq:sym3mn}) is the identity endomorphism
(the last summand requires $n\geq 3$ to be nonzero). In particular, all three projections are nonzero when $n\geq 3$.

For $n\geq2$, the following shows that in any {symmetric rank} decomposition of~$\SMn$, it is impossible to have all  summands corresponding to matrices $L_i$ of rank one.  Moreover, for $n\geq3$, at least one summand
corresponds to a matrix having rank at least $3$.  
We note that this statement is in contrast to 
tensor decompositions of~$\SMn$ where
there do exist decompositions constructed from rank one tensors.  In fact, the \mbox{$n=3$} case is in stark 
contrast to usual matrix multiplication
where there exist decompositions 
for which no matrix appearing has rank greater than 
one, e.g., {the standard decomposition}.

\begin{thm0}\label{thm:ki}
Suppose that $\SMn=\sum_{i=1}^k  L_i^3$
is a symmetric rank decomposition.
If $n=2$, there exists $i$ with $ \rank(L_i) = 2$.
Moreover, if $n\ge 3$, $\max_i \rank(L_i)\ge 3$. 
\end{thm0}
\begin{proof}
Any  summand $L_i^3$ with $\trank L_i=1$ is of  the form $L_i=v_i\otimes \omega_i\in V\otimes V^\vee$
and induces an element of rank one  
that takes $a\otimes b\otimes c$ to  $\omega_i(a)\omega_i(b)\omega_i(c)v_i^3$
which vanishes outside $\mathrm{Sym}^3V$.
This element lies in  $\tend(\mathrm{Sym}^3V)$ in the decomposition~\eqref{eq:sym3mn}. Hence, any sum of these elements  lies in this subspace and thus projects to zero in the second and third factors
in (\ref{eq:sym3mn}).

Similarly, any summand of rank two only gives rise to a term appearing in 
$\mathrm{Sym}^3V\ot \mathrm{Sym}^3V^*\op S_{21}V\ot S_{21}V^* $ because
one needs three independent vectors for a term in $\La 3 V\ot \La 3 V^*$.
\end{proof}

This following provides a 
slight improvement over the na\"\i ve bound
of $8n^3$.

\begin{prop0}[A modest upper bound]
${\bold R_s}(\SMn)\le 8{n\choose 3}+4{n\choose 2}+n$.
\end{prop0}
\begin{proof}
Every monomial appearing in $\SMn$ has the form
$a_{ij}a_{jk}a_{ki}$.  This bound arises
from considering the symmetric ranks of each of these monomials.
There are~$2{n\choose 3}$ monomials corresponding to distinct cardinality 3 sets
$\{i,j,k\}\subset\{1,\ldots, n\}$ and each monomial has symmetric rank $4$. 
There are $2{n\choose 2}$ monomials corresponding to distinct cardinality 2 sets
$\{i,j\}\subset\{1,\ldots, n\}$
and they group together in ${n\choose 2}$ pairs
as $a_{ij}a_{ji}(a_{ii}+a_{jj})$ with each such term having symmetric rank two. 
Finally, there are $n$ monomials of the form $a_{ii}^3$ for $i=1,\ldots, n$.
\end{proof}

The following considers algebraic geometric aspects
of $\SMn$.
 
\begin{prop0}\label{prop:smnirred} 
(i) The singular locus of $\{ \SMn=0\} \subset \BP \Matn$ is
$$\{ [A]\in \BP \Matn \mid A^2=0\}.$$

(ii)  The polynomial $\SMtwo$ is reducible, while $\SMn$ is irreducible for $n\ge 3$.

\end{prop0}
\begin{proof}
The map $(A,B)\mapsto\mathrm{tr}(AB^t)$ is a nondegenerate pairing.
Since we can write $\mathrm{tr}(A^3)=\mathrm{tr}(A\cdot A^2)$,
this proves $(i)$. 

Note that the $(i,j)$ entry of $A^2$ coincides, up to scalar multiple, with the partial derivative $\dfrac{\partial\left(\SMn\right)}{\partial a_{j,i}}$. 
 In order to prove $(ii)$, we estimate the dimension of the singular locus computed in $(i)$.
If $A$ belongs to the singular locus  of $\{ \SMn=0\}$,
we know $\ker(A)\subseteq \mathrm{im}(A)$
so that $\rank(A)\le n/2$.
It follows that the singular locus of $\{ \SMn=0\}$ has codimension $\ge 3$ for $n\ge 3$ 
showing that $\SMn$ must be irreducible.
If not, the singular locus contains the intersection of any two irreducible components, having
codimension $\le 2$. The $n=2$ case follows from (\ref{eq:f2}) below.
\end{proof}

\subsection{Decomposition of $\SMtwo$}\label{subsec:SMtwo} 

The reducibility of $\SMtwo$ is as follows:

\begin{eqnarray}
\SMtwo =& {a}_{0,0}^{3}+3 {a}_{0,0} {a}_{0,1} {a}_{1,0}+3 {a}_{0,1}
      {a}_{1,0} {a}_{1,1}+{a}_{1,1}^{3} \nonumber\\
=&\underbrace{\ttrace(A)}_{\textrm{non tg hyperp.}}\cdot
\underbrace{\left(\ttrace^2(A)-3\det(A)\right)}_{\textrm{smooth quadric}}.\label{eq:f2}
\end{eqnarray}
In particular, for this classically studied polynomial,
its zero set is the union of 
a smooth quadric and a non-tangent hyperplane. 
A general cubic surface has a unique Waring decomposition as a sum of $5$ summands by the 
Sylvester Pentahedral Theorem \cite[Theor. 3.9]{OeOtt}.
Hence, every $f\in\mathrm{Sym}^3\C^4$ has $\ur_s(f )\le 5$.
However, $\bold R_s(\SMtwo)=6$ (see \cite[IV, \S 97]{Seg}) with a minimal Waring decomposition given by 
\begin{equation}\label{eq:6f2} \SMtwo= \sum_{i=1}^6  L_i^3\end{equation}
where
\vskip 0.3cm
$L_1=\dfrac 12 \begin{pmatrix}-1&1\\-1&-1\end{pmatrix}$,
$L_2=\dfrac 12\begin{pmatrix}-1&-1\\1&-1\end{pmatrix}$,
$L_3=\dfrac 12\begin{pmatrix}1&1\\1&1\end{pmatrix}$,
$L_4=\dfrac 12\begin{pmatrix}1&-1\\-1&1\end{pmatrix}$,
\vskip 0.3cm
$L_5= \begin{pmatrix}1&0\\0&0\end{pmatrix}$, and 
$L_6= \begin{pmatrix}0&0\\0&1\end{pmatrix}$.

\begin{rem0}[A remark on 5 summands]
For the decomposition presented in \eqref{eq:6f2},
$\trank (L_i)=2$ for $i=1, 2$
while $\trank (L_i)=1$ for $i=3,4,5,6$
in agreement with Theorem~\ref{thm:ki}.
Since $\SMtwo$ is $GL_2$-invariant for the {\it conjugate action} which takes
$A$ to $G^{-1}AG$ for every $G\in GL_2$, 
the matrices $L_i$ can be replaced in \eqref{eq:6f2} with $G^{-1} L_iG$ for any $G\in GL_2$.

Consider a family $f_{{2},\epsilon}$ 
which has a Waring decomposition given by five {matrices $L_{i,\epsilon}$} for $\epsilon\neq 0$
and $f_{{2},0}=sM_{\langle {2}\rangle}$.
In all the examples we have found,
the five  {matrices $L_{i,\epsilon}$} converge 
as $\epsilon\to 0$ to the identity matrix
that is indeed a fixed point for the conjugate action.
\end{rem0}

The following Remark provides a geometric description 
for decompositions of $\SMtwo$ using six terms. 

\begin{rem0}
Identify  the projective space of $2\times 2$ matrices with $\mathbb P^3$. Let $\Qc$ 
be the quadric of matrices of rank $1$ 
and let $\ell$ denote the line spanned by the identity $I$ and the skew-symmetric point $\Lambda$.

For a choice of $3$ points $Q_1,Q_2,Q_3$ in the intersection of $\Qc$ with the plane of traceless matrices,
let $A_1,A_2,B_1,B_2,C_1,C_2$ denote the 6 points  of intersection of the two rulings of $\Qc$
passing through each $Q_i$.  
These points, together with $I$, determine a minimal decomposition
of the general tensor $\Mtwo$, as explained in \cite{CILO}.

A decomposition of $\SMtwo$ is determined as follows: let $Q_3$ be the intersection of the lines
$(B_1C_1)$ and $(B_2C_2)$. Then the six points $L_1\ldots L_6$ are obtained by taking
$L_6=A_2$, $L_5=A_1$, $L_4=$ the intersection of $(B_1,C_1)$ with the plane $\pi$ of symmetric matrices, 
$L_3=$ the intersection of $(B_2,C_2)$ with $\pi$, $L_2=$ the intersection of 
the line $(Q_3A_2)$ with $\ell$ (they meet), and $L_1=$ the intersection of the line
$(Q_3,A_1)$ with $\ell$.

For instance, starting with
$Q_1= \begin{pmatrix} 0&1\\ 0& 0\end{pmatrix}$,
$Q_2=\begin{pmatrix}  0& 0\\1& 0\end{pmatrix}$, and
$Q_3=\begin{pmatrix} 1&1\\-1&-1\end{pmatrix}$,
we obtain the six points $L_1,\dots,L_6$ of the decomposition \eqref{eq:6f2} described above. 
\end{rem0}

We ask if an analogous geometric description could provide small decompositions of $\SMn$ for $n\geq3$.

\subsection{Case of $\SMthree$}
The polynomial $\SMthree$
is irreducible by Proposition \ref{prop:smnirred} with the following lower bound on border rank.

\begin{prop0}\label{prop:14}
$\ur_s(\SMthree)\ge 14.$
\end{prop0}
\begin{proof}
Let $V=\C^9$. 
For any $\phi\in \mathrm{Sym}^3V$ we have the linear map 
{$$A_{\phi}\colon V^\vee\otimes\Lambda^4V\to\Lambda^5V\otimes V$$} which is defined by contracting the
elements of the source with $\phi$ and then projecting to the target.
This projection is well-defined because the map 
{$$S^2V\ot \La 4 V\rightarrow \La 5 V\ot V$$ }
is a 
$GL(V)$-module map and the image of the projection is
the unique copy of  {$S_{411}V\subset \La 5 V\ot V$}. 
This map was denoted as {$\mathrm{YF}_{3,8}(\phi)$}  in \cite[Eq.~(2)]{LOann}.
Using \cite{GS}, direct computation shows
that $\trank A_{v^3}=70$ for a nonzero $v \in V$ and 
$\trank A_{\SMthree}=950$.
By linearity
$\ur_s(\SMthree)\ge \lceil\frac{950}{70}\rceil = 14.$
\end{proof}

The following provides information
on the rank. 

\begin{thm0*}\label{thm*:ms3_19}
${\bold R}_s(\SMthree) \leq 18$
with a Waring decomposition of 
$\SMthree$ with $18$ summands 
found numerically by
{\tt Bertini} \cite{BertiniSoftware}
with all $18$ summations having rank~$3$.
\end{thm0*}
\begin{proof}
After numerically approximating a decomposition
with {\tt Bertini} \cite{BertiniSoftware},
applying the isosingular local dimension test \cite{Isosingular} 
suggested that there is at least one $9$-dimensional family of decompositions.
We used the extra $9$ degrees of freedom
to set $9$ entries to $0$, $1$, or $-1$
producing a polynomial system which has an isolated
nonsingular root with an approximation
given in Appendix~\ref{appendix}
and electronically available at 
\url{http://dx.doi.org/10.7274/R0VT1Q1J}.
\end{proof}

{Decompositions with $18$ summands
were highly structured
leading to the following.}

\begin{conj}\label{conj}
${\bold R}_s(\SMthree) = 18$.
\end{conj}

In our experiments, 
we were unable to compute a decomposition
of~$\SMthree$ using $18$ summands 
with real matrices.

\section{The polynomials $\SMsn$ and $\SMson$}
 \label{sec:smsn}

We start with statements from Theorem~\ref{thm:main}.

\begin{proof}[Proof of \eqref{SMsnom} and \eqref{SMsoom}]
The following two inequalities are
trivial since $\SMsn$ is a specialization of $\SMn$:
$${\bold R}_s(\SMsn)\le{\bold R}_s(\SMn)
\hbox{~~~and~~~} {\bold R}(\SMsn)\le{\bold R}(\SMn).$$

For $n\times n$ matrices $A$, $B$, $C$ consider the $3n \times 3n$ symmetric matrix 
$$X =\begin{pmatrix}0& C^T& A\\
 C& 0& B^T\\
 A^T& B& 0\end{pmatrix}.$$
We have $\ttrace(X^3) = 6\ttrace(ABC)$ since
$$X^3=\begin{pmatrix}ABC+C^TB^TA^T& *& *\\
 *& CAB+B^TA^TC^T& *\\
 *& *& BCA+A^TC^TB^T\end{pmatrix}.$$
It immediately follows ${\bold R}(\Mn)\le{\bold R}(sM^S_{\langle 3n\rangle})$.  Hence, \eqref{SMsnom} follows by
  a similar argument as in the proof of \eqref{SMnom}.

Since $X$ is traceless, the same argument also proves  \eqref{SMsoom}. 
\end{proof}

\def\nothing{
The following is analogous to Theorem~\ref{thm:ki}. 

\begin{thm0}\label{thm:ski}
Suppose that $\SMsn=\sum_{i=1}^k  L_i^3$
is a symmetric rank decomposition.
If $n=2$, there exists $i$ with $ \rank(L_i) = 2$.
Moreover, if $n\ge 3$, $\max_i \rank(L_i)\ge 3$. 
\end{thm0}
\begin{proof}
We have
$$\mathrm{Sym}^3\mathrm{Sym}^2V =\mathrm{Sym}^6V\oplus S_{4,2}V\oplus S_{2,2,2}V.
$$ 
The element $\SMsn$ projects to the identity
in each of the three spaces following
the same approach as for $\SMn$.
We omit the rest of the proof since it is analogous to the proof of the Theorem~\ref{thm:ki}.
\end{proof}
}

\subsection{Decomposition of $\sSMtwo$}\label{subsec:smstwo}

As in the general case \eqref{eq:f2}, $\sSMtwo$ is a reducible polynomial while $\SMsn$ is irreducible  for $n\ge 3$
(the same argument as in Proposition \ref{prop:smnirred} works).
In fact, 
$$\sSMtwo\begin{pmatrix}a_0&a_1\\a_1&a_2\end{pmatrix}=(a_0+a_2)(a_0^2+3a_1^2-a_0a_2+a_2^2),$$ 
which corresponds to the union of a smooth conic with a secant (not tangent) line. 
Moreover, it was known classically
that $\ur_s(\sSMtwo) = \bold R_s(\sSMtwo)=4$, which 
  is the generic rank in $\PP(\mathrm{Sym}^3\C^3)$
  with a minimal Waring decomposition given~by 
\begin{equation}\label{eq:6f22}
6\cdot \sSMtwo=L_1^3+L_2^3-2L_3^3-2L_4^3\end{equation}
where
$$
{\tiny
L_1=\begin{pmatrix}2&\frac{\sqrt{-1}}{2}\\\frac{\sqrt{-1}}{2}&0\end{pmatrix},
L_2=\begin{pmatrix}0&-\frac{\sqrt{-1}}{2}\\-\frac{\sqrt{-1}}{2}&2\end{pmatrix},
L_3=\begin{pmatrix}1&\sqrt{-1}\\\sqrt{-1}&0\end{pmatrix},
L_4=\begin{pmatrix}0&-\sqrt{-1}\\-\sqrt{-1}&1\end{pmatrix}}.$$
We note that $L_1$ and $L_2$ are similar 
as well as $L_3$ and $L_4$
all have rank $2$.
 
\subsection{Case of $sM^{S}_{\langle 3\rangle}$}\label{subsec:tableau} 

{We consider $\ur_s(sM^{S}_{\langle 3\rangle})$ as a cubic polynomial on $\C^6$.}
Since the generic rank in $\PP(\mathrm{Sym}^3\C^6)$ 
is $10$ (see \cite{AH1995}),
we have $\ur_s(sM^{S}_{\langle 3\rangle})\leq 10$.
To show that equality holds,  consider
the degree $10$ invariant in $\mathrm{Sym}^{10}(\mathrm{Sym}^3\C^{6})$ corresponding to the following 
Young diagram (see, e.g., \cite[\S 3.9]{5inv}, for the symbolic notation of invariants):
  \begin{equation}\label{eq:t10}T_{10}=\begin{matrix}\young({1}1122,23334,44555,66767,78898,9\dieci 9\dieci\dieci).\end{matrix}\end{equation}
This invariant is generalized in \cite[Prop.~3.25]{BI:2017}.

\begin{lemma0}\label{lem:t10} The polynomial $T_{10}$ defined by \eqref{eq:t10} forms 
a basis of the $SL_6$-invariant space $(\mathrm{Sym}^{10}(\mathrm{Sym}^3\C^{6}))^{SL_6}$ and
is in the ideal of $\s_9(\nu_3(\pp 5))$. 
Moreover, $T_{10}(sM^{S}_{\langle 3\rangle})\neq 0$
showing that $\ur_s(sM^{S}_{\langle 3\rangle}) > 9$.
\end{lemma0}

\begin{proof} A plethysm calculation, e.g.,
using \texttt{Schur} \cite{schur},
shows that 
$$\dim (\mathrm{Sym}^{10}(\mathrm{Sym}^3\C^{6}))^{SL_6}=1.$$
We explicitly evaluated $T_{10}(sM^{S}_{\langle 3\rangle})$ using the same algorithm as in \cite{16051} and \cite{CIM:17}
which phrases the evaluation as a tensor contraction and   ignores summands that contribute zero to the result.
The result was that $T_{10}(sM^{S}_{\langle 3\rangle})\neq0$.

We now consider evaluating $T_{10}$ on all cubics 
of the form $f=\sum_{i=1}^9\ell_i^3$.  
The expression $T_{10}(f)$ splits as the sum of several
terms of the form
$T_{10}(\ell_{i_1}^3,\ldots,\ell_{i_{10}}^3)$ where, in each of these summands, there is a repetition of 
some $\ell_i$.
We claim that every $T_{10}(\ell_{i_1}^3,\ldots,\ell_{i_{10}}^3)$ vanishes due to this repetition. Indeed, each pair $(i,j)$ with $1\le i<j\le 10$
appears in at least one column of \eqref{eq:t10}. In other words, for any $g\colon \{1,\ldots,10\}\to \{1,\ldots,9\}$, the tableau evaluation $g(T_{10})$
has a repetition in at least one column and thus vanishes. This approach is the main tool used in \cite{16051}.
\end{proof}

Since the polynomial $T_{10}$
vanishes on $\sigma_9(\nu_3(\BP^5))$,
Lemma \ref{lem:t10} immediately
yields that $\ur_s(sM^{S}_{\langle 3\rangle}) = 10$, because $\sigma_{10}(v_3(\BP^5))$ equals
$\BP \mathrm{Sym}^3\C^{6}$.
The following considers decompositions
with $10$ summands.

\begin{prop0*}\label{prop*:10}
${\bold R}_s(sM^{S}_{\langle 3\rangle}) = 10$.
 \end{prop0*}
\begin{proof}
Consider decompositions consisting
of 3 symmetric matrices each of the form
\begin{equation}\label{eq:10rank2}
{\small
\begin{pmatrix}
* & * & * \\
* & 0 & 0 \\
* & 0 & 0 
\end{pmatrix},~~
\begin{pmatrix}
0 & * & 0\\
* & * & * \\
0 & * & 0 
\end{pmatrix},~~
\begin{pmatrix}
0 & 0 & *\\
0 & 0 & * \\
* & * & * 
\end{pmatrix},
}
\end{equation}
each of which is clearly rank deficient, and one 
symmetric matrix of the form
\begin{equation}\label{eq:10traceless}
{\small
\begin{pmatrix}
0 & * & *\\
* & 0 & * \\
* & * & 0 
\end{pmatrix}
}
\end{equation}
which is clearly traceless.  

Upon substituting these forms, which have a total
of $3\cdot 10 = 30$ unknowns, 
into the $\binom{5+3}{3}=56$ equations which describe
the decompositions of $sM^{S}_{\langle 3\rangle}$,
there are $28$ equations which vanish identically
leaving $28$ polynomial equations in $30$ affine
variables.
The isosingular local dimension test \cite{Isosingular}
in {\tt Bertini} \cite{BertiniSoftware}
suggests that this system 
has at least one $3$-dimensional solution component
which we utilize the~$3$ extra degrees of freedom
to make one entry either $\pm 1$ 
in one of each of the three types of matrices in \eqref{eq:10rank2}.  The resulting system has an isolated solution 
which we present one here to $4$ significant digits:
$$
{\tiny
\begin{array}{l}
\begin{pmatrix}0.1755 & 2.16 & -1\\2.16 & 0 & 0\\-1 & 0 & 0\\ \end{pmatrix},
\begin{pmatrix}0.6889 & -0.4607 & -0.8745\\-0.4607 & 0& 0\\-0.8745 & 0& 0\\\end{pmatrix},
\begin{pmatrix}0.874 & 0.1991 & 0.5836\\0.1991 & 0& 0\\0.5836 & 0& 0\\\end{pmatrix},\\[0.3cm]
\begin{pmatrix}0 & -0.7877 & 0\\-0.7877 & 0.5269 & 1\\0 & 1 & 0\\\end{pmatrix},
\begin{pmatrix}0 & 1.431 & 0\\1.431 & 0.326 & 0.9555\\0 & 0.9555 & 0\\\end{pmatrix},
\begin{pmatrix}0 & 0.076 & 0\\0.076 & 0.9356 & -0.4331\\0 & -0.4331 & 0\\\end{pmatrix},\\[0.3cm]
\begin{pmatrix}0 & 0& 1\\0 & 0& 0.2278\\1 & 0.2278 & 0.6677\\\end{pmatrix},
\begin{pmatrix}0& 0& -0.6362\\0& 0& 0.4255\\-0.6362 & 0.4255 & 0.8077\\\end{pmatrix},
\begin{pmatrix}0 & 0& -0.09825\\0& 0& -1.21\\-0.09825 & -1.21 & 0.5599\\\end{pmatrix},\\[0.3cm]
\begin{pmatrix}0 & -2.317 & 0.8998\\-2.317 & 0& -0.4797\\0.8998 & -0.4797 & 0\\\end{pmatrix}.
\end{array}
}
$$
The eigenvalues $\lambda_1,\lambda_2,\lambda_3$ of the 
first $9$ summands satisfy
$$\lambda_3 = (\lambda_1 + \lambda_2)(\lambda_1^2 + \lambda_1 \lambda_2 + \lambda_2^2) - 1 = 0$$
while the eigenvalues of the traceless matrix satisfy
$$\lambda_1+\lambda_2+\lambda_3 = 
(\lambda_1+\lambda_2)(\lambda_1+\lambda_3)(\lambda_2+\lambda_3)+2 = 0.$$
\end{proof}

The variety $\s_9(\nu_3(\pp 5))$ has codimension $2$
as expected. The following describes
generators of its ideal.

\begin{thm0*}\label{thm:9secant}
The variety $\s_9(\nu_3(\pp 5))$
has codimension $2$
and degree $280$.  It is the 
complete intersection of the solution set of 
$T_{10}$ and a hypersurface of degree $28$.
\end{thm0*}
\begin{proof}
It is easy to computationally verify that 
the variety $X := \sigma_{9}(\nu_3(\BP^5))\subset\BP^{55}$
has the expected dimension of $53$, e.g., via \cite[Lemma~3]{Projections}.  
This also follows from the Alexander-Hirschowitz Theorem \cite{AH1995}. We used the approach
in \cite[\S 2]{LowerBound} with 
{\tt Bertini}~\cite{BertiniSoftware} 
to compute a so-called pseudowitness set \cite{Projections}
for $X$ yielding \mbox{$\deg X = 280$}.
With this pseudowitness set, \cite{aCM,aG}
shows that $X$ is arithmetically Cohen-Macaulay
and arithmetically Gorenstein.
In particular, the Hilbert function of 
the finite set $X\cap \sL$ where $\sL\subset\BP^{55}$ is a general linear space of dimension $2$ is
$$
\begin{array}{r} 1,3,6,10,15,21,28,36,45,55,65,75,85,95,
105,115,125,135,145,155,165,175,\\
~~~~~~185,195,205,215,225,235,244,
252,259,265,270,274,277,279,280,280.
\end{array}
$$
Thus, the ideal of $X\cap\sL$ is minimally
generated by a degree $10$ polynomial (corresponding to $T_{10}$)
and a polynomial of degree $28$.  
The same holds for $X$, i.e., $X$ is a complete
intersection defined by the vanishing of $T_{10}$
and a polynomial of degree~$28$,
since $X$ is arithmetically Cohen-Macaulay.
The Hilbert series of $X$ is
$$
\frac{\hbox{$
\begin{array}{r}
1 + 2 t + 3 t^{2} + 4 t^{3} + 5 t^{4} + 6 t^{5} + 7 t^{6} + 8 t^{7} + 9 t^{8} + 10 t^{9} (1 + t + t^2 + \cdots + t^{18})
\\
~~~+~9 t^{28} + 8 t^{29} + 7 t^{30} + 6 t^{31} + 5 t^{32} + 4 t^{33} + 3 t^{34} + 2 t^{35} + t^{36} \end{array}$}}{(1-t)^{54}}
$$
\end{proof}
 
\begin{rem0}  The generic polynomial in $\s_9(\nu_3(\pp 5))$ has exactly two Waring decompositions,
this is the last subgeneric  $\s_k(\nu_d(\pp n))$ whose generic member has a  non-unique Waring decomposition \cite[Thm. 1.1]{COV2}.
\end{rem0}

We close with the traceless $3\times 3$ case
$sM^{S,0}_{\langle 3\rangle}$ where we take
$a_5 = -(a_0+a_3)$.  

\begin{prop0}\label{prop:8}
$\ur_s(sM^{S,0}_{\langle 3\rangle})= 
\bold R_s(sM^{S,0}_{\langle 3\rangle}) = 8$
\end{prop0}
\begin{proof}
Although $\sigma_7(\nu_3(\BP^4))\subset\BP^{34}$
is expected to fill the ambient space, it  
is defective: it is 
a hypersurface of degree $15$ defined by
the cubic root of the determinant of a 45$\times$45 matrix, e.g., see \cite{Ott:09,16051}.
This $45\times45$ matrix 
evaluated at $sM^{S,0}_{\langle 3\rangle}$
has full rank showing that
$\ur_s(sM^{S,0}_{\langle 3\rangle}) > 7$. 
Since $8$ is the generic rank, $\ur_s(sM^{S,0}_{\langle 3\rangle}) = 8$.  

To show the existence of a decomposition using
$8$ summands, we need to solve
a system of $\binom{4+3}{3} = 35$ polynomials
in $40$ affine variables.  By including the 
determinant of the matrices corresponding
to the first $5$ summands,
we produce a square system with $40$ polynomials
in $40$ variables.  We prove the existence of a solution
via \mbox{$\alpha$-theory} using {\tt alphaCertified} \cite{alphaCertified} starting with the following 
approximation:
$$
{\tiny
\begin{array}{l}
\begin{pmatrix} 0.2533609-0.3253227i & 0.3900781-0.4785431i & 0.2123864-0.4078949i\\0.3900781-0.4785431i & 2.017554-0.09428536i & 0.2851566+1.393898i\\0.2123864-0.4078949i & 0.2851566+1.393898i & -2.2709149+0.41960806i\\\end{pmatrix},\\[0.3cm]
\begin{pmatrix} 0.04310556+0.1656553i & 0.1274312+0.3981205i & -0.4116999-0.01601833i\\0.1274312+0.3981205i & -1.934801+0.2834521i & -0.4336855-1.882461i\\-0.4116999-0.01601833i & -0.4336855-1.882461i & 1.89169544-0.4491074i\\\end{pmatrix},\\[0.3cm]
\begin{pmatrix} -0.3785169+0.4133459i & 0.250335-1.167081i & 0.2879607-0.03026006i\\0.250335-1.167081i & 0.3783613-0.2548455i & -0.2320977-0.1928574i\\0.2879607-0.03026006i & -0.2320977-0.1928574i & 0.00015556-0.1585004i\\\end{pmatrix},\\[0.3cm]
\begin{pmatrix} 0.3088717+0.01475953i & -0.1783686-0.3906188i & 0.2316816+0.3868949i\\-0.1783686-0.3906188i & 2.09271-0.2084929i & 0.4260931+1.812166i\\0.2316816+0.3868949i & 0.4260931+1.812166i & -2.4015817+0.19373337i\\\end{pmatrix},\\[0.3cm]
\begin{pmatrix} 0.1338705-0.583525i & 0.8378364-1.288966i & 0.9626108+0.2184638i\\0.8378364-1.288966i & 0.3123518-0.5020508i & 0.9180592+0.8550508i\\0.9626108+0.2184638i & 0.9180592+0.8550508i & -0.4462223+1.0855758i\\\end{pmatrix},\\[0.3cm]
\begin{pmatrix} 0.1972212-0.8229062i & -0.1848761-0.6480481i & 0.1550351-0.03293642i\\-0.1848761-0.6480481i & 1.063376-1.941261i & 1.409034+0.4606352i\\0.1550351-0.03293642i & 1.409034+0.4606352i & -1.2605972+2.7641672i\\\end{pmatrix},\\[0.3cm]
\begin{pmatrix} 0.8841674-0.3162212i & 1.394277+0.2241248i & 0.1553027+0.6211393i\\1.394277+0.2241248i & 0.7687433-0.1439827i & -0.1014701+1.35298i\\0.1553027+0.6211393i & -0.1014701+1.35298i & -1.6529107+0.4602039i\\\end{pmatrix},\\[0.3cm]
\begin{pmatrix} 0.1094107+0.06367402i & -0.3608902+1.394814i & -0.7766249-0.2283304i\\-0.3608902+1.394814i & -0.5622866+0.5667869i & -0.01617272+0.05988273i\\-0.7766249-0.2283304i & -0.01617272+0.05988273i & 0.4528759-0.63046092i\\\end{pmatrix}
\end{array}
}
$$
where $i = \sqrt{-1}$.
\end{proof}

\section{The polynomial $\ZMn$}
\label{sec:szmn}
Let $Z_n$ be the space of symmetric matrices with 
zeros on the diagonal. 
The cubic $\ZMn(A)$ is a   polynomial in ${{n}\choose 2}$ indeterminates,
its na\"\i ve expression has ${{n}\choose 3}$ terms: 
\begin{equation}\label{eq:binomialn3}\ZMn(A)=\sum_{1\le i<j<k\le n}a_{ij}a_{jk}a_{ik}\end{equation}
The proof of \eqref{ZMnom} is similar to the others
and thus omitted.

Since $sM_{\langle 2 \rangle}^{Z}$ is identically zero,
we take $n\geq 3$.  
Let $P_n$ denote the finite set of~$2^{n-1}$ vectors 
of the form 
\mbox{$v=\left(1,\pm 1,\ldots, \pm 1\right)^T\in\Z^n$}.
In Theorem \ref{thm:waringzmn}, we use~$P_n$ 
to construct a decomposition of $\ZMn$.
Although such a decomposition
is not minimal for $n\ge 6$ (see Proposition~\ref{thm:waringzmn2}), a modification of it 
constructs the decomposition~\eqref{eq:waringzm8} 
for $n=8$ which we believe is minimal (see Remark \ref{rem:waringzmn3}).

For each $v\in P_n$, $vv^T-I_n\in Z_n$
with eigenvalues $\{-1,\ldots, -1,n-1\}$
and off-diagonal elements $\pm 1$.

\begin{thm0}\label{thm:waringzmn} For $n\ge 3$, 
we have the decomposition of $2^{n-1}$ summands:
\begin{equation}\label{eq:waringzmn}
2^{n+2}\ZMn(A) = \sum_{v\in P_n}\left(\ttrace[(vv^T-I_n)\cdot A]\right)^3.
\end{equation}
\end{thm0}
\begin{proof} If $v=\left(v_1,\ldots, v_n\right)^T\in P_n$, 
then the $(i,j)$-entry of $(vv^T-I_n)$ is $v_iv_j-\delta_{ij}$.
The monomials appearing in $\left(\ttrace[(vv^T-I_n)\cdot A]\right)^3=\left(\sum_{i<j}a_{ij}v_iv_j\right)^3$ divide into three groups
\begin{enumerate}
\item $a_{ij}^3v_iv_j$ for $i<j$,
\item $a_{ij}^2a_{pq}v_pv_q$ for $i<j$ and $p<q$,
\item $a_{ij}a_{pq}a_{rs}v_iv_jv_pv_qv_rv_s$.
\end{enumerate}

Summing over $P_n$, the monomials of the first group cancel each other because, for any fixed value $v_i\in\{-1, +1\}$,  the vectors $v\in P_n$
having this fixed value divide into two subsets of equal size, having respectively $v_j=-1$ or $v_j=1$. This argument includes the case $i=0$,
when $v_0=1$.

For the same reason the monomials of the second group cancel each other.

In the third group, all monomials when $\#\{i,j,p,q,r,s\}\ge 4$ cancel each other because there is an index which appear only once,
and the above argument shows that the sum over this index makes zero.
If   $\#\{i,j,p,q,r,s\}=3$ and the monomial is not in the first or second group, then each index appears exactly twice and we get exactly all the summands which appear in (\ref{eq:binomialn3}).

Since these cover all cases, the right-hand side of \eqref{eq:waringzmn} sums up to a scalar multiple
of the left-hand side.
\end{proof}

\begin{prop0}\label{thm:waringzmn1}
The decompositions in (\ref{eq:waringzmn}) are minimal for $n=3,4,5$.
In particular, $\ur_s(sM_{\langle n \rangle}^{Z})
= \bold R_s(sM_{\langle n \rangle}^{Z}) = 2^{n-1}$ for $n=3, 4, 5$.
\end{prop0}

\begin{proof}
We compute the  Koszul flattening $\mathrm{YF}_{3,n^2-1}(sM_{\langle n \rangle}^{Z})$ as in \cite[(2)]{LOann}
where~$r_n$ is its rank. Let $q_n=\rank\mathrm{YF}_{3,n^2-1}(\ell^3)={{m
_n}\choose{\lfloor m_n/2\rfloor}}$ for any linear form $\ell$, where $m_n=n(n-1)/2-1$.
With this setup, \cite[Prop. 4.1.1]{LOann}
and Theorem~\ref{thm:waringzmn}
provide
$$\left\lceil\frac{r_n}{q_n}\right\rceil \leq \ur_s(sM_{\langle n \rangle}^{Z}) \leq
\bold R_s(sM_{\langle n \rangle}^{Z}) \leq 2^{n-1}.$$
The result follows immediately from the following table:
$$\begin{array}{r|r|r|r|r}n&r_n&q_n&\left\lceil\frac{r_n}{q_n}\right\rceil & 2^{n-1}\\
\hline
3 &    8 &   2 &  4 & 4\\
4 &   72 &  10 &  8 & 8\\
5 & 1920 & 126 & 16 & 16\end{array}$$
\end{proof}

For comparison, the known lower bounds
on the border rank of $\Mn$ when $n = 3,4,5$
are $15,29,47$, respectively, 
with the general lower bound from~\cite{2016arXiv160807486L}
of $\ur(\Mn)\geq 2n^2-\tlog_2(n)-1$.

For $n=6,7$, the decomposition \eqref{eq:waringzmn} has $32$ and $64$ summands, respectively. 
The following shows that such decompositions
are not minimal (as well for any $n\ge 6$). 
\begin{prop0}\label{thm:waringzmn2} $\bold R_s(sM_{\langle 6 \rangle}^{Z}) \le 30$ and
$\bold R_s(sM_{\langle 7 \rangle}^{Z}) \le 48$.
\end{prop0}
\begin{proof} 
See Appendix~\ref{appendix}
for a decomposition of $32sM_{\langle 6 \rangle}^{Z}$ with $30$ summands having integer coefficients 
and a decomposition of $160sM_{\langle 7 \rangle}^{Z}$ with $48$ summands having coefficients in $\mathbb Q[\sqrt{5}]$.
\end{proof}

When $n = 8$, the following provides
a decomposition using $64$ summands.

\begin{prop0}\label{pro:waringzmn3} 
Let $P_8^+$ be the subset of $P_8$ consisting of $v$ such that $+1$ appears an even number of times, so 
$\# P_8^+ = 64$.
Then, 
\begin{equation}\label{eq:waringzm8}
2^9 sM_{\langle 8 \rangle}^{Z}(A) = \sum_{v\in P_8^+}\left(\ttrace[(vv^T-I_n)\cdot A]\right)^3\end{equation}
\end{prop0}
\begin{proof} This is easy to verify by direct computation.
\end{proof}

\begin{rem0}\label{rem:waringzmn3} 
We expect the decomposition \eqref{eq:waringzm8} 
with $64$ summands is minimal.
In the $\binom{8}{2} = 28$ indeterminants
of the matrix $A$, consider
the polynomial $f(A;\ell_1,\ldots,\ell_{64})= \sum_{i=1}^{64}\left(\ttrace[(\ell_i^T A]\right)^3$.
Evaluated at 
$\ell_i=(v_iv_i^T-I_n)$ where \mbox{$v_i\in P_8^+$},
the polynomial $f$ has maximal rank of $1792 = 28\cdot64$.
\end{rem0}

Proposition~\ref{pro:waringzmn3} suggests
one should look for strategic subsets 
$P_n^{?}\subset P_n$ analogous to $P_8^+\subset P_8$
to produce minimal decompositions.  
For $n=9$ and $10$, we can use $P_9^+$ and $P_{10}^+$
to obtain (again by direct computation)
$\bold R_s(sM_{\langle 9 \rangle}^{Z}) \le 128$ 
and $\bold R_s(sM_{\langle 10 \rangle}^{Z}) \le 256$, 
but both seem to not be sharp.

 \appendix
 
 \section{Decompositions}\label{appendix}

The following $18$ matrices 
of rank $3$ form a numerical approximation of a decomposition of $\SMthree$ where $i = \sqrt{-1}$:
$${\tiny
\begin{array}{l}
\begin{pmatrix} -0.13-0.311i & 0.499-0.51i & -0.464-0.387i\\-1.4-2.08i & 2.46-0.687i & -1.56+0.414i\\-0.141-0.542i & 0.44-0.374i & -0.783-0.0408i\end{pmatrix},\\[0.3cm]
\begin{pmatrix} 0.568+1.31i & -0.592+0.375i & 1\\0.0129-0.785i & 0.598+0.73i & -1.48+0.0928i\\0.943-0.486i & 0.407+0.64i & -0.55-0.572i\end{pmatrix},\\[0.3cm]
\begin{pmatrix} -0.557-0.103i & 0.169-0.756i & 0.198+0.804i\\0.815-1.25i & 1 & -1\\1.23+0.517i & -0.491+0.197i & 0.516-1.16i\end{pmatrix},\\[0.3cm]
\begin{pmatrix} -0.649-0.377i & 0.787+0.21i & 0\\1.26-1.57i & 1.2+2.02i & -0.712-1.79i\\-0.314+0.107i & 0.602+0.0423i & 0.0664-0.178i\end{pmatrix},\\[0.3cm]
\begin{pmatrix} 0.714+0.0554i & 0.283-0.0242i & -0.0436-1.28i\\0.491+2.16i & -0.449-0.276i & 2.3-1.63i\\0.685+1.21i & -0.692-0.311i & 0.695-1.04i\end{pmatrix},\\[0.3cm]
\begin{pmatrix} -1.34-0.753i & -0.344-0.339i & -0.0879+1.74i\\0.00563-2.43i & -0.0178-0.303i & -1.91+1.15i\\-0.148-0.755i & 0.106+0.39i & -0.312+0.239i\end{pmatrix},\\
\end{array}
}$$
$${\tiny
\begin{array}{l}
\begin{pmatrix} -1.42-0.99i & 0.779-0.573i & -1.33+0.129i\\-1.36-0.599i & 0.496-0.462i & -0.474-0.937i\\-1.5-6.93\cdot 10^{-5}i & 0.209-0.712i & -0.747+0.635i\end{pmatrix},\\[0.3cm]
\begin{pmatrix} 0.918+0.932i & -0.867-0.478i & 1+0.368i\\-1.01+0.753i & 0.132-1.7i & -0.195+2.31i\\1& -0.659+0.171i & 0.493-0.273i\end{pmatrix},\\[0.3cm]
\begin{pmatrix} 1.22+1.7i & -0.408+0.328i & 0.739-1.56i\\-0.731+2.33i & -0.0271+0.0704i & 1.36-0.778i\\0.137+0.165i & 0.395-0.0697i & -0.575-0.311i\end{pmatrix},\\[0.3cm]
\begin{pmatrix} 1.5+0.508i & 0.406+0.256i & 0.672-0.572i\\0.665-0.0681i & 1.5+0.508i & -0.529+0.896i\\0.224+0.652i & 0.000481-0.584i & 1.02+0.317i\end{pmatrix},\\[0.3cm]
\begin{pmatrix} 0.0701-0.426i & 0.128+0.459i & 0.0339+0.782i\\-0.144-1.56i & 0.409+0.605i & -1.92+1.94i\\-1.12-0.862i & 0.596-0.152i & -0.351+1.68i\end{pmatrix},\\[0.3cm]
\begin{pmatrix} -1.25+1.2i & -0.0161-1.08i & 0.231+1.64i\\-0.909+0.146i & 0.96-0.305i & -1.86+0.264i\\-0.134+1.44i & -0.622-0.495i & 0.903+0.567i\end{pmatrix},\\[0.3cm]
\begin{pmatrix} -1 & -0.09-0.772i & -1.17+0.839i\\-1.95+0.717i & -0.092-1.51i & -0.511+0.706i\\-0.905+0.354i & -0.208-0.553i & -0.484+1.31i\end{pmatrix},\\[0.3cm]
\begin{pmatrix} -1.44-1.13i & 0.414-0.261i & -0.828+1.13i\\1& -0.751+0.559i & 0.0548-1.41i\\0.0382+0.716i & -0.481+0.194i & 0.519-0.242i\end{pmatrix},\\[0.3cm]
\begin{pmatrix} 0 & 0.0696+0.285i & 0.537+0.0341i\\0.321+0.252i & 0 & -0.612+0.169i\\-0.178+0.381i & 0.248-0.256i & -0.126-0.284i\end{pmatrix},\\[0.3cm]
\begin{pmatrix} 0.261-0.0237i & -0.258+0.311i & -0.338-0.577i\\-1.21-0.256i & -0.221-0.63i & -0.347+0.516i\\-0.129-1.47i & 0.789+0.348i & -1.71-0.163i\end{pmatrix},\\[0.3cm]
\begin{pmatrix} -0.294+1.33i & -0.471-0.0831i & 0.353-1.03i\\-0.411+1.21i & -0.534-0.0483i & 1.34-1.84i\\1.27+0.0632i & -0.0116+0.723i & -0.748-1.48i\end{pmatrix},\\[0.3cm]
\begin{pmatrix} -1.47+1.27i & -0.701-0.533i & 1.66+0.192i\\-1.88+0.176i & -0.351-0.513i & 1.59+0.427i\\0.568+0.624i & -0.306+0.631i & 0.153-1.58i\end{pmatrix}.
\\
\end{array}
}$$

A decomposition of $32 sM^Z_{\langle 6\rangle}$ using
$30$ summands is:
$$
{\tiny
\begin{array}{l}
(a_{2} - 2 a_{3} + a_{4} - a_{6} - a_{8} + 2 a_{9} - a_{10} - a_{12} - a_{13} - a_{15})^3 + \\
(2 a_{1} + a_{2} - 2 a_{3} - a_{4} + a_{6} - 2 a_{7} - a_{8} - a_{10} - 2 a_{11} - a_{12} + a_{13} + a_{15})^3 + \\
(2 a_{3} - a_{2} - a_{4} - a_{6} + 2 a_{7} - a_{8} - a_{10} + 2 a_{11} - a_{12} - a_{13} + 2 a_{14} - a_{15})^3 + \\
(a_{2} - 2 a_{1} + 2 a_{3} - a_{4} - a_{6} + a_{8} - 2 a_{9} + a_{10} + a_{12} - a_{13} + 2 a_{14} - a_{15})^3 + \\
(2 a_{1} + a_{2} - 2 a_{3} + a_{4} + 2 a_{5} + a_{6} - 2 a_{7} + a_{8} + 2 a_{9} - a_{10} + a_{12} - a_{13} - 2 a_{14} + a_{15})^3 + \\
(2 a_{1} - a_{2} - a_{4} - a_{6} - a_{8} - a_{10} + a_{12} - a_{13} - 2 a_{14} + a_{15})^3 + \\
(a_{2} + a_{4} - 2 a_{5} + a_{6} + a_{8} - 2 a_{9} - a_{10} + 2 a_{11} - a_{12} - a_{13} + 2 a_{14} - a_{15})^3 + \\
(a_{2} + 2 a_{3} - a_{4} - 2 a_{5} + a_{6} + 2 a_{7} - a_{8} - 2 a_{9} + a_{10} - a_{12} - a_{13} + a_{15})^3 + \\
(2 a_{3} - a_{2} + a_{4} - 2 a_{5} - a_{6} + 2 a_{7} + a_{8} - 2 a_{9} - a_{10} + a_{12} + a_{13} - a_{15})^3 + \\
(a_{4} - a_{2} - 2 a_{1} - 2 a_{5} + a_{6} + 2 a_{7} - a_{8} + a_{10} + a_{12} - a_{13} + 2 a_{14} - a_{15})^3 + \\
(2 a_{1} + a_{2} + a_{4} + a_{6} + a_{8} + a_{10} - a_{12} + a_{13} - 2 a_{14} - a_{15})^3 + \\
(a_{2} - a_{4} - a_{6} - 2 a_{7} + a_{8} + 2 a_{9} + a_{10} - 2 a_{11} - a_{12} - a_{13} - 2 a_{14} + a_{15})^3 + \\
(a_{2} + a_{4} + 2 a_{5} - a_{6} - 2 a_{7} - a_{8} + a_{10} + a_{12} + a_{13} + a_{15})^3 + \\
(2 a_{1} - a_{2} - 2 a_{3} + a_{4} - a_{6} - 2 a_{7} + a_{8} + a_{10} - 2 a_{11} + a_{12} - a_{13} - a_{15})^3 + \\
(2 a_{3} - a_{2} - 2 a_{1} - a_{4} - 2 a_{5} + a_{6} + a_{8} - a_{10} + 2 a_{11} + a_{12} - a_{13} + a_{15})^3 + \\
(a_{6} - 2 a_{3} - a_{4} - a_{2} + a_{8} + 2 a_{9} + a_{10} + a_{12} + a_{13} + a_{15})^3 + \\
(a_{6} - a_{2} - a_{4} - 2 a_{1} + 2 a_{7} + a_{8} - 2 a_{9} + a_{10} + 2 a_{11} - a_{12} + a_{13} - a_{15})^3 + \\
(a_{4} - a_{2} + a_{6} - 2 a_{7} - a_{8} + 2 a_{9} - a_{10} - 2 a_{11} + a_{12} + a_{13} - 2 a_{14} - a_{15})^3 + \\
(2 a_{5} - a_{4} - a_{2} + a_{6} - 2 a_{7} + a_{8} - a_{10} - a_{12} - a_{13} - a_{15})^3 + \\
(2 a_{3} - a_{2} - 2 a_{1} + a_{4} + a_{6} - a_{8} - 2 a_{9} - a_{10} - a_{12} + a_{13} + 2 a_{14} + a_{15})^3 + \\
(2 a_{1} + a_{2} - a_{4} + 2 a_{5} + a_{6} - a_{8} + 2 a_{9} + a_{10} - 2 a_{11} + a_{12} - a_{13} - a_{15})^3 + \\
(2 a_{1} - a_{2} + a_{4} + 2 a_{5} - a_{6} + a_{8} + 2 a_{9} - a_{10} - 2 a_{11} - a_{12} + a_{13} + a_{15})^3 + \\
(a_{2} - 2 a_{1} + a_{4} - a_{6} + 2 a_{7} - a_{8} - 2 a_{9} - a_{10} + 2 a_{11} + a_{12} - a_{13} + a_{15})^3 + \\
(a_{10} - a_{4} - 2 a_{5} - a_{6} - a_{8} - 2 a_{9} - a_{2} + 2 a_{11} + a_{12} + a_{13} + 2 a_{14} + a_{15})^3 + \\
(2 a_{1} - a_{2} - 2 a_{3} - a_{4} + 2 a_{5} - a_{6} - 2 a_{7} - a_{8} + 2 a_{9} + a_{10} - a_{12} + a_{13} - 2 a_{14} - a_{15})^3 + \\
(a_{4} - 2 a_{3} - a_{2} + 2 a_{5} + a_{6} - a_{8} + a_{10} - 2 a_{11} - a_{12} - a_{13} - 2 a_{14} + a_{15})^3 + \\
(a_{2} + 2 a_{3} + a_{4} + a_{6} + 2 a_{7} + a_{8} + a_{10} + 2 a_{11} + a_{12} + a_{13} + 2 a_{14} + a_{15})^3 + \\
(a_{2} - 2 a_{1} - a_{4} - 2 a_{5} - a_{6} + 2 a_{7} + a_{8} - a_{10} - a_{12} + a_{13} + 2 a_{14} + a_{15})^3 + \\
(a_{2} - 2 a_{1} + 2 a_{3} + a_{4} - 2 a_{5} - a_{6} - a_{8} + a_{10} + 2 a_{11} - a_{12} + a_{13} - a_{15})^3 + \\
(a_{2} - 2 a_{3} - a_{4} + 2 a_{5} - a_{6} + a_{8} - a_{10} - 2 a_{11} + a_{12} + a_{13} - 2 a_{14} - a_{15})^3.
\end{array}
}
$$

With $\beta = \sqrt{5}/2$, a decomposition of $160 sM^Z_{\langle 7\rangle}$ using
$48$ summands is:
$$
{\tiny
\begin{array}{l}
\begin{array}{l}
(a_{3} - a_{1} - 2 a_{4} + a_{5} - a_{6} - a_{8} + 2 a_{9} - a_{10} + a_{11} - 2 a_{16} + a_{17} - a_{18} - 2 a_{19} + 2 a_{20} - a_{21}~~~~~ \\[-0.01cm]
~~~-~2 a_{2} \beta + 2 a_{7} \beta - 2 a_{12} \beta + a_{13} \beta - 2 a_{14} \beta + 2 a_{15} \beta)^3 + \end{array}\\[-0.01cm]
\begin{array}{l}
(2 a_{3} - 2 a_{1} + a_{4} + 2 a_{5} - 2 a_{6} - 2 a_{8} - a_{9} - 2 a_{10} + 2 a_{11} + a_{16} + 2 a_{17} - 2 a_{18} + a_{19} - a_{20} - 2 a_{21})^3 + \end{array} \\[-0.01cm]
\begin{array}{l}
(a_{4} - 2 a_{3} - 2 a_{1} - 2 a_{5} - 2 a_{6} + 2 a_{8} - a_{9} + 2 a_{10} + 2 a_{11} - a_{16} + 2 a_{17} + 2 a_{18} - a_{19} - a_{20} + 2 a_{21})^3 + \end{array} \\[-0.01cm]
\begin{array}{l}
(a_{1} - a_{3} + 2 a_{4} + a_{5} + a_{6} - a_{8} + 2 a_{9} + a_{10} + a_{11} - 2 a_{16} - a_{17} - a_{18} + 2 a_{19} + 2 a_{20} + a_{21} \\[-0.01cm]
~~~+~2 a_{2} \beta + 2 a_{7} \beta - 2 a_{12} \beta + a_{13} \beta + 2 a_{14} \beta + 2 a_{15} \beta)^3 + \end{array} \\[-0.01cm]
\begin{array}{l}
(2 a_{1} + 2 a_{3} - a_{4} - 2 a_{5} + 2 a_{6} + 2 a_{8} - a_{9} - 2 a_{10} + 2 a_{11} - a_{16} - 2 a_{17} + 2 a_{18} + a_{19} - a_{20} - 2 a_{21})^3 + \end{array} \\[-0.01cm]
\begin{array}{l}
(2 a_{3} - 2 a_{1} - a_{4} - 2 a_{5} - 2 a_{6} - 2 a_{8} + a_{9} + 2 a_{10} + 2 a_{11} - a_{16} - 2 a_{17} - 2 a_{18} + a_{19} + a_{20} + 2 a_{21})^3 + \end{array} \\[-0.01cm]
\begin{array}{l}
(a_{1} - a_{3} - 2 a_{4} - a_{5} + a_{6} - a_{8} - 2 a_{9} - a_{10} + a_{11} + 2 a_{16} + a_{17} - a_{18} + 2 a_{19} - 2 a_{20} - a_{21} \\[-0.01cm]
~~~+~2 a_{2} \beta + 2 a_{7} \beta - 2 a_{12} \beta - a_{13} \beta - 2 a_{14} \beta + 2 a_{15} \beta)^3 + \end{array} \\[-0.01cm]
\begin{array}{l}
(a_{1} - a_{3} + 2 a_{4} - a_{5} - a_{6} - a_{8} + 2 a_{9} - a_{10} - a_{11} - 2 a_{16} + a_{17} + a_{18} - 2 a_{19} - 2 a_{20} + a_{21} \\[-0.01cm]
~~~+~2 a_{2} \beta + 2 a_{7} \beta - 2 a_{12} \beta + a_{13} \beta - 2 a_{14} \beta - 2 a_{15} \beta)^3 + \end{array} \\[-0.01cm]
\begin{array}{l}
(a_{1} + a_{3} - 2 a_{4} - a_{5} - a_{6} + a_{8} - 2 a_{9} - a_{10} - a_{11} - 2 a_{16} - a_{17} - a_{18} + 2 a_{19} + 2 a_{20} + a_{21} \\[-0.01cm]
~~~+~ 2 a_{2} \beta + 2 a_{7} \beta + 2 a_{12} \beta - a_{13} \beta - 2 a_{14} \beta - 2 a_{15} \beta)^3 + \end{array} \\[-0.01cm]
\begin{array}{l}
(2 a_{4} - a_{3} - a_{1} + a_{5} - a_{6} + a_{8} - 2 a_{9} - a_{10} + a_{11} - 2 a_{16} - a_{17} + a_{18} + 2 a_{19} - 2 a_{20} - a_{21} \\[-0.01cm]
~~~+~2 a_{2} \beta - 2 a_{7} \beta - 2 a_{12} \beta + a_{13} \beta + 2 a_{14} \beta - 2 a_{15} \beta)^3 + \end{array} \\[-0.01cm]
\begin{array}{l}
(a_{1} + a_{3} - 2 a_{4} + a_{5} + a_{6} + a_{8} - 2 a_{9} + a_{10} + a_{11} - 2 a_{16} + a_{17} + a_{18} - 2 a_{19} - 2 a_{20} + a_{21} \\[-0.01cm]
~~~+~2 a_{2} \beta + 2 a_{7} \beta + 2 a_{12} \beta - a_{13} \beta + 2 a_{14} \beta + 2 a_{15} \beta)^3 + \end{array} \\[-0.01cm]
\begin{array}{l}
(a_{3} - a_{1} + 2 a_{4} - a_{5} - a_{6} - a_{8} - 2 a_{9} + a_{10} + a_{11} + 2 a_{16} - a_{17} - a_{18} - 2 a_{19} - 2 a_{20} + a_{21} \\[-0.01cm]
~~~-~2 a_{2} \beta + 2 a_{7} \beta - 2 a_{12} \beta - a_{13} \beta + 2 a_{14} \beta + 2 a_{15} \beta)^3 + \end{array} \\[-0.01cm]
\begin{array}{l}
(a_{1} + a_{3} + 2 a_{4} + a_{5} - a_{6} + a_{8} + 2 a_{9} + a_{10} - a_{11} + 2 a_{16} + a_{17} - a_{18} + 2 a_{19} - 2 a_{20} - a_{21} \\[-0.01cm]
~~~-~ 2 a_{2} \beta - 2 a_{7} \beta - 2 a_{12} \beta - a_{13} \beta - 2 a_{14} \beta + 2 a_{15} \beta)^3 + \end{array} \\[-0.01cm]
\begin{array}{l}
(a_{3} - a_{1} + 2 a_{4} - a_{5} - a_{6} - a_{8} - 2 a_{9} + a_{10} + a_{11} + 2 a_{16} - a_{17} - a_{18} - 2 a_{19} - 2 a_{20} + a_{21} \\[-0.01cm]
~~~+~ 2 a_{2} \beta - 2 a_{7} \beta + 2 a_{12} \beta + a_{13} \beta - 2 a_{14} \beta - 2 a_{15} \beta)^3 + \end{array} \\[-0.01cm]
\begin{array}{l}
(2 a_{4} - a_{3} - a_{1} - a_{5} + a_{6} + a_{8} - 2 a_{9} + a_{10} - a_{11} - 2 a_{16} + a_{17} - a_{18} - 2 a_{19} + 2 a_{20} - a_{21}  \\[-0.01cm]
~~~-~2 a_{2} \beta + 2 a_{7} \beta + 2 a_{12} \beta - a_{13} \beta + 2 a_{14} \beta - 2 a_{15} \beta)^3 + \end{array} \\[-0.01cm]
\begin{array}{l}
(a_{1} + a_{3} + 2 a_{4} - a_{5} + a_{6} + a_{8} + 2 a_{9} - a_{10} + a_{11} + 2 a_{16} - a_{17} + a_{18} - 2 a_{19} + 2 a_{20} - a_{21} \\[-0.01cm]
~~~-~ 2 a_{2} \beta - 2 a_{7} \beta - 2 a_{12} \beta - a_{13} \beta + 2 a_{14} \beta - 2 a_{15} \beta)^3 + \end{array} \\[-0.01cm]
\begin{array}{l}
(a_{8} - a_{3} - 2 a_{4} - a_{5} - a_{6} - a_{1} + 2 a_{9} + a_{10} + a_{11} + 2 a_{16} + a_{17} + a_{18} + 2 a_{19} + 2 a_{20} + a_{21} \\[-0.01cm]
~~~+~ 2 a_{2} \beta - 2 a_{7} \beta - 2 a_{12} \beta - a_{13} \beta - 2 a_{14} \beta - 2 a_{15} \beta)^3 + \end{array} \\[-0.01cm]
\begin{array}{l}
(2 a_{3} - 2 a_{1} + a_{4} - 2 a_{5} + 2 a_{6} - 2 a_{8} - a_{9} + 2 a_{10} - 2 a_{11} + a_{16} - 2 a_{17} + 2 a_{18} - a_{19} + a_{20} - 2 a_{21})^3 + \end{array} \\[-0.01cm]
\begin{array}{l}
(a_{3} - a_{1} - 2 a_{4} - a_{5} + a_{6} - a_{8} + 2 a_{9} + a_{10} - a_{11} - 2 a_{16} - a_{17} + a_{18} + 2 a_{19} - 2 a_{20} - a_{21} \\[-0.01cm]
~~~+~ 2 a_{2} \beta - 2 a_{7} \beta + 2 a_{12} \beta - a_{13} \beta - 2 a_{14} \beta + 2 a_{15} \beta)^3 + \end{array}\\[-0.01cm]
\begin{array}{l}
(2 a_{1} - 2 a_{3} - a_{4} + 2 a_{5} + 2 a_{6} - 2 a_{8} - a_{9} + 2 a_{10} + 2 a_{11} + a_{16} - 2 a_{17} - 2 a_{18} - a_{19} - a_{20} + 2 a_{21})^3 + \end{array} \\[-0.01cm]
\begin{array}{l}
(2 a_{4} - a_{3} - a_{1} - a_{5} + a_{6} + a_{8} - 2 a_{9} + a_{10} - a_{11} - 2 a_{16} + a_{17} - a_{18} - 2 a_{19} + 2 a_{20} - a_{21} \\[-0.01cm]
~~~+~ 2 a_{2} \beta - 2 a_{7} \beta - 2 a_{12} \beta + a_{13} \beta - 2 a_{14} \beta + 2 a_{15} \beta)^3 + \end{array} \\[-0.01cm]
\begin{array}{l}
(a_{3} - a_{1} + 2 a_{4} + a_{5} + a_{6} - a_{8} - 2 a_{9} - a_{10} - a_{11} + 2 a_{16} + a_{17} + a_{18} + 2 a_{19} + 2 a_{20} + a_{21} \\[-0.01cm]
~~~-~ 2 a_{2} \beta + 2 a_{7} \beta - 2 a_{12} \beta - a_{13} \beta - 2 a_{14} \beta - 2 a_{15} \beta)^3 + \end{array} \\[-0.01cm]
\begin{array}{l}
(a_{3} - a_{1} + 2 a_{4} + a_{5} + a_{6} - a_{8} - 2 a_{9} - a_{10} - a_{11} + 2 a_{16} + a_{17} + a_{18} + 2 a_{19} + 2 a_{20} + a_{21} \\[-0.01cm]
~~~+~ 2 a_{2} \beta - 2 a_{7} \beta + 2 a_{12} \beta + a_{13} \beta + 2 a_{14} \beta + 2 a_{15} \beta)^3 + \end{array} \\[-0.01cm]
\begin{array}{l}
(2 a_{3} - 2 a_{1} - a_{4} + 2 a_{5} + 2 a_{6} - 2 a_{8} + a_{9} - 2 a_{10} - 2 a_{11} - a_{16} + 2 a_{17} + 2 a_{18} - a_{19} - a_{20} + 2 a_{21})^3 + \end{array} \\[-0.01cm]
\begin{array}{l}
(a_{3} - a_{1} - 2 a_{4} + a_{5} - a_{6} - a_{8} + 2 a_{9} - a_{10} + a_{11} - 2 a_{16} + a_{17} - a_{18} - 2 a_{19} + 2 a_{20} - a_{21} \\[-0.01cm]
~~~+~ 2 a_{2} \beta - 2 a_{7} \beta + 2 a_{12} \beta - a_{13} \beta + 2 a_{14} \beta - 2 a_{15} \beta)^3 + \end{array} \\[-0.01cm]
\begin{array}{l}
(a_{1} + a_{3} - 2 a_{4} + a_{5} + a_{6} + a_{8} - 2 a_{9} + a_{10} + a_{11} - 2 a_{16} + a_{17} + a_{18} - 2 a_{19} - 2 a_{20} + a_{21} \\[-0.01cm]
~~~-~ 2 a_{2} \beta - 2 a_{7} \beta - 2 a_{12} \beta + a_{13} \beta - 2 a_{14} \beta - 2 a_{15} \beta)^3 + \end{array} \\[-0.01cm]
\begin{array}{l}
(a_{3} - a_{1} - 2 a_{4} - a_{5} + a_{6} - a_{8} + 2 a_{9} + a_{10} - a_{11} - 2 a_{16} - a_{17} + a_{18} + 2 a_{19} - 2 a_{20} - a_{21} \\[-0.01cm]
~~~-~ 2 a_{2} \beta + 2 a_{7} \beta - 2 a_{12} \beta + a_{13} \beta + 2 a_{14} \beta - 2 a_{15} \beta)^3 + \end{array} \\[-0.01cm]
\begin{array}{l}
(2 a_{1} - 2 a_{3} + a_{4} + 2 a_{5} - 2 a_{6} - 2 a_{8} + a_{9} + 2 a_{10} - 2 a_{11} - a_{16} - 2 a_{17} + 2 a_{18} + a_{19} - a_{20} - 2 a_{21})^3 + \end{array} \\[-0.01cm]
\begin{array}{l}
(a_{1} + a_{3} - 2 a_{4} - a_{5} - a_{6} + a_{8} - 2 a_{9} - a_{10} - a_{11} - 2 a_{16} - a_{17} - a_{18} + 2 a_{19} + 2 a_{20} + a_{21} \\[-0.01cm]
~~~-~ 2 a_{2} \beta - 2 a_{7} \beta - 2 a_{12} \beta + a_{13} \beta + 2 a_{14} \beta + 2 a_{15} \beta)^3 + \end{array} \\[-0.01cm]
\begin{array}{l}
(a_{1} - a_{3} + 2 a_{4} - a_{5} - a_{6} - a_{8} + 2 a_{9} - a_{10} - a_{11} - 2 a_{16} + a_{17} + a_{18} - 2 a_{19} - 2 a_{20} + a_{21} \\[-0.01cm]
~~~-~ 2 a_{2} \beta - 2 a_{7} \beta + 2 a_{12} \beta - a_{13} \beta + 2 a_{14} \beta + 2 a_{15} \beta)^3 + \end{array} \\[-0.01cm]
\begin{array}{l}
(2 a_{4} - a_{3} - a_{1} + a_{5} - a_{6} + a_{8} - 2 a_{9} - a_{10} + a_{11} - 2 a_{16} - a_{17} + a_{18} + 2 a_{19} - 2 a_{20} - a_{21} \\[-0.01cm]
~~~-~ 2 a_{2} \beta + 2 a_{7} \beta + 2 a_{12} \beta - a_{13} \beta - 2 a_{14} \beta + 2 a_{15} \beta)^3 + \end{array} \\[-0.01cm]
\begin{array}{l}
(a_{5} - a_{3} - 2 a_{4} - a_{1} + a_{6} + a_{8} + 2 a_{9} - a_{10} - a_{11} + 2 a_{16} - a_{17} - a_{18} - 2 a_{19} - 2 a_{20} + a_{21} \\[-0.01cm]
~~~-~ 2 a_{2} \beta + 2 a_{7} \beta + 2 a_{12} \beta + a_{13} \beta - 2 a_{14} \beta - 2 a_{15} \beta)^3 + \end{array} \\[-0.01cm]
\begin{array}{l}
(a_{1} - a_{3} - 2 a_{4} + a_{5} - a_{6} - a_{8} - 2 a_{9} + a_{10} - a_{11} + 2 a_{16} - a_{17} + a_{18} - 2 a_{19} + 2 a_{20} - a_{21} \\[-0.01cm]
~~~-~2 a_{2} \beta - 2 a_{7} \beta + 2 a_{12} \beta + a_{13} \beta - 2 a_{14} \beta + 2 a_{15} \beta)^3 + \end{array} \\[-0.01cm]
\begin{array}{l}
(a_{1} + a_{3} + 2 a_{4} - a_{5} + a_{6} + a_{8} + 2 a_{9} - a_{10} + a_{11} + 2 a_{16} - a_{17} + a_{18} - 2 a_{19} + 2 a_{20} - a_{21} \\[-0.01cm]
~~~+~ 2 a_{2} \beta + 2 a_{7} \beta + 2 a_{12} \beta + a_{13} \beta - 2 a_{14} \beta + 2 a_{15} \beta)^3 + \end{array} \\[-0.01cm]
\begin{array}{l}
(2 a_{5} - 2 a_{3} - a_{4} - 2 a_{1} - 2 a_{6} + 2 a_{8} + a_{9} - 2 a_{10} + 2 a_{11} + a_{16} - 2 a_{17} + 2 a_{18} - a_{19} + a_{20} - 2 a_{21})^3 + \end{array} \\[-0.01cm]
\begin{array}{l}
(a_{8} - a_{3} - 2 a_{4} - a_{5} - a_{6} - a_{1} + 2 a_{9} + a_{10} + a_{11} + 2 a_{16} + a_{17} + a_{18} + 2 a_{19} + 2 a_{20} + a_{21} \\[-0.01cm]
~~~-~ 2 a_{2} \beta + 2 a_{7} \beta + 2 a_{12} \beta + a_{13} \beta + 2 a_{14} \beta + 2 a_{15} \beta)^3 + \end{array} \\[-0.01cm]
\begin{array}{l}
(2 a_{1} + 2 a_{3} + a_{4} + 2 a_{5} + 2 a_{6} + 2 a_{8} + a_{9} + 2 a_{10} + 2 a_{11} + a_{16} + 2 a_{17} + 2 a_{18} + a_{19} + a_{20} + 2 a_{21})^3 + \end{array} \\[-0.01cm]
\begin{array}{l}
(a_{5} - a_{3} - 2 a_{4} - a_{1} + a_{6} + a_{8} + 2 a_{9} - a_{10} - a_{11} + 2 a_{16} - a_{17} - a_{18} - 2 a_{19} - 2 a_{20} + a_{21} \\[-0.01cm]
~~~+~ 2 a_{2} \beta - 2 a_{7} \beta - 2 a_{12} \beta - a_{13} \beta + 2 a_{14} \beta + 2 a_{15} \beta)^3 + \end{array} \\[-0.01cm]
\begin{array}{l}
(2 a_{1} + 2 a_{3} - a_{4} + 2 a_{5} - 2 a_{6} + 2 a_{8} - a_{9} + 2 a_{10} - 2 a_{11} - a_{16} + 2 a_{17} - 2 a_{18} - a_{19} + a_{20} - 2 a_{21})^3 + \end{array} \\[-0.01cm]
\begin{array}{l}
(2 a_{6} - 2 a_{3} - a_{4} - 2 a_{5} - 2 a_{1} + 2 a_{8} + a_{9} + 2 a_{10} - 2 a_{11} + a_{16} + 2 a_{17} - 2 a_{18} + a_{19} - a_{20} - 2 a_{21})^3 + \end{array} \\[-0.01cm]
\begin{array}{l}
(2 a_{1} + 2 a_{3} + a_{4} - 2 a_{5} - 2 a_{6} + 2 a_{8} + a_{9} - 2 a_{10} - 2 a_{11} + a_{16} - 2 a_{17} - 2 a_{18} - a_{19} - a_{20} + 2 a_{21})^3 + \end{array} \\[-0.01cm]
\begin{array}{l}
(a_{1} - a_{3} + 2 a_{4} + a_{5} + a_{6} - a_{8} + 2 a_{9} + a_{10} + a_{11} - 2 a_{16} - a_{17} - a_{18} + 2 a_{19} + 2 a_{20} + a_{21} \\[-0.01cm]
~~~-~ 2 a_{2} \beta - 2 a_{7} \beta + 2 a_{12} \beta - a_{13} \beta - 2 a_{14} \beta - 2 a_{15} \beta)^3 + \end{array}\\[-0.01cm]
\begin{array}{l}
(2 a_{1} - 2 a_{3} - a_{4} - 2 a_{5} - 2 a_{6} - 2 a_{8} - a_{9} - 2 a_{10} - 2 a_{11} + a_{16} + 2 a_{17} + 2 a_{18} + a_{19} + a_{20} + 2 a_{21})^3 + \end{array} \\[-0.01cm]
\begin{array}{l}
(a_{4} - 2 a_{3} - 2 a_{1} + 2 a_{5} + 2 a_{6} + 2 a_{8} - a_{9} - 2 a_{10} - 2 a_{11} - a_{16} - 2 a_{17} - 2 a_{18} + a_{19} + a_{20} + 2 a_{21})^3 + \end{array} \\[-0.01cm]
\begin{array}{l}
(a_{1} + a_{3} + 2 a_{4} + a_{5} - a_{6} + a_{8} + 2 a_{9} + a_{10} - a_{11} + 2 a_{16} + a_{17} - a_{18} + 2 a_{19} - 2 a_{20} - a_{21} \\[-0.01cm]
~~~+~ 2 a_{2} \beta + 2 a_{7} \beta + 2 a_{12} \beta + a_{13} \beta + 2 a_{14} \beta - 2 a_{15} \beta)^3 + \end{array}\\[-0.01cm]
\begin{array}{l}
(2 a_{1} - 2 a_{3} + a_{4} - 2 a_{5} + 2 a_{6} - 2 a_{8} + a_{9} - 2 a_{10} + 2 a_{11} - a_{16} + 2 a_{17} - 2 a_{18} - a_{19} + a_{20} - 2 a_{21})^3 + \end{array} \\[-0.01cm]
\begin{array}{l}
(a_{1} - a_{3} - 2 a_{4} + a_{5} - a_{6} - a_{8} - 2 a_{9} + a_{10} - a_{11} + 2 a_{16} - a_{17} + a_{18} - 2 a_{19} + 2 a_{20} - a_{21} \\[-0.01cm]
~~~+~ 2 a_{2} \beta + 2 a_{7} \beta - 2 a_{12} \beta - a_{13} \beta + 2 a_{14} \beta - 2 a_{15} \beta)^3. \end{array} 
\end{array}
}
$$
 
\end{document}